\documentclass[12pt,leqno]{article}
\usepackage{amsmath,amssymb,amscd,latexsym}
\newcommand{\C}{{\mathbb{C}}}
\newcommand{\F}{{\mathbb{F}}}
\newcommand{\Q}{{\mathbb{Q}}}
\newcommand{\R}{{\mathbb{R}}}
\newcommand{\Sa}{{\mathbb{S}}}
\newcommand{\Z}{{\mathbb{Z}}}
\newcommand{\Co}{\mathrm{Co}}
\newcommand{\et}{\mathrm{\acute{e}t}}
\newcommand{\id}{\mathrm{id}}
\newcommand{\imm}{\mathrm{im}\,}
\newcommand{\ind}{\mathrm{ind}}
\renewcommand{\mod}{\;\mathrm{mod}\;}
\newcommand{\prim}{\mathrm{prim}}
\newcommand{\proj}{\mathrm{proj}}
\newcommand{\sgn}{\mathrm{sgn}\,}
\newcommand{\spec}{\mathrm{spec}\,}
\newcommand{\vol}{\mathrm{vol}}
\newcommand{\Frob}{\mathrm{Frob}}
\newcommand{\GL}{\mathrm{GL}\,}
\newcommand{\Imm}{\mathrm{Im}\,}
\newcommand{\Ker}{\mathrm{Ker}\,}
\newcommand{\PSO}{\mathrm{PSO}}
\newcommand{\RRe}{\mathrm{Re}\,}
\newcommand{\SO}{\mathrm{SO}}
\renewcommand{\Sp}{\mathrm{Sp}}
\newcommand{\Tr}{\mathrm{Tr}}
\newcommand{\tr}{\mathrm{tr}}
\newcommand{\Ah}{{\mathcal A}}
\newcommand{\Dh}{{\mathcal D}}
\newcommand{\Fh}{{\mathcal F}}
\newcommand{\Hh}{{\mathcal H}}
\newcommand{\Lh}{{\mathcal L}}
\newcommand{\Oh}{{\mathcal O}}
\newcommand{\Rh}{{\mathcal R}}
\newcommand{\eo}{\mathfrak{o}}
\newcommand{\ep}{\mathfrak{p}}
\newcommand{\eX}{{\mathcal X}}
\newcommand{\oH}{\bar{H}}
\newcommand{\ohne}{\setminus}
\newcommand{\silo}{\stackrel{\sim}{\longrightarrow}}
\newcommand{\tei}{\, | \,}
\newcommand{\ent}{\;\widehat{=}\;}
\newcommand{\hullet}{\raisebox{0.05cm}{$\scriptscriptstyle \bullet$}}
\newcommand{\verk}{\mbox{\scriptsize $\,\circ\,$}}
\newcommand{\halb}{\frac{1}{2}}

\newtheorem{theorem}{Theorem}[section]
\newtheorem{cor}[theorem]{Corollary}
\newenvironment{example}{{\bf Example}}{}
\newtheorem{remarks}[theorem]{Remarks}
\newenvironment{remarknn}{\noindent {\bf Remark}}{}
\newtheorem{conj}[theorem]{Conjecture}
\newtheorem{fact}[theorem]{Fact}
\newtheorem{punkt}[theorem]{$\!\!$}
\newenvironment{proofof}{\noindent {\bf Proof of}}{\mbox{}\hfill$\Box$}
\begin{document}
\title{Number theory and dynamical systems on foliated spaces}
\author{Christopher Deninger}
\date{\ }
\maketitle
\section{Introduction}
In this paper we report on some developments in the search for a dynamical understanding of number theoretical zeta functions that have taken place since my ICM lecture \cite{D2}. We also point out a number of problems in analysis that will have to be solved in order to make further progress.

In section 2 we give a short introduction to foliations and their cohomology. Section 3 is devoted to progress on the dynamical Lefschetz trace formula for one-codimensional foliations mainly due to \'Alvarez L\'opez and Kordyukov. In section 4 we make the comparison with the ``explicit formulas'' in analytic number theory. Finally in section 5 we generalize the conjectural dynamical Lefschetz trace formula of section 3 to phase spaces which are more general than manifolds. This was suggested by the number theoretical analogies of section 4. 

This account is written from an elementary point of view as far as arithmetic geometry is concerned, in particular motives are not mentioned. In spirit the present article is therefore a sequel to \cite{D1}.

There is a different approach to number theoretical zeta functions using dynamical systems by A. Connes \cite{Co}. His phase space is a non-commutative quotient of the ad\`eles. Although superficially related, the two approaches seem to be deeply different. Whereas Connes' approach generalizes readily to automorphic $L$-functions \cite{So} but not to motivic $L$-functions, it is exactly the opposite with our picture. One may wonder whether there is some kind of Langlands correspondence between the two approaches.

I would like to thank the Belgium and German mathematical societies very much for the opportunity to lecture about this material during the joint BMS--DMV meeting in Li\`ege 2001.
\section{Foliations and their cohomology}
A $d$-dimensional foliation $\Fh = \Fh_X$ on a smooth manifold $X$ of dimension $a$ is a partition of $X$ into immersed connected $d$-dimensional manifolds $F$, the ``leaves''. Locally the induced partition should be trivial: Every point of $X$ should have an open neighborhood $U$ diffeomorphic to an open ball $B$ in $\R^a$ such that the leaves of the induced partition on $U$ correspond to the submanifolds $B \cap (\R^d \times \{ y \})$ of $B$ for $y$ in $\R^{a-d}$.

One of the simplest non-trivial examples is the one-dimensional foliation of the two-dimensional torus $T^2 = \R^2 / \Z^2$ by lines of irrational slope $\alpha$. These are given by the immersions
\[
\R \hookrightarrow T^2 \; , \; t \mapsto (x + t \alpha , t) \mod \Z^2
\]
parametrized by $x \mod \Z + \alpha \Z$. In this case every leaf is dense in $T^2$ and the intersection of a global leaf with a small open neighborhood $U$ as above decomposes into countably many connected components. It is the global behaviour which makes foliations complicated. For a comprehensive introduction to foliation theory, the reader may turn to \cite{Go} for example.

To a foliation $\Fh$ on $X$ we may attach its tangent bundle $T \Fh$ whose total space is the union of the tangent spaces to the leaves. By local triviality of the foliation it is a sub vector bundle of the tangent bundle $TX$. It is integrable i.e. the commutator of any two vector fields with values in $T \Fh$ again takes values in $T\Fh$. Conversely a theorem of Frobenius asserts that every integrable sub vector bundle of $TX$ arises in this way.

Differential forms of order $n$ along the leaves are defined as the smooth sections of the real vector bundle $\Lambda^n T^* \Fh$, 
\[
\Ah^n_{\Fh} (X) = \Gamma (X, \Lambda^n T^* \Fh) \; .
\]
The same formulas as in the classical case define exterior derivatives along the leaves:
\[
d^n_{\Fh} : \Ah^n_{\Fh} (X) \longrightarrow \Ah^{n+1}_{\Fh} (X) \; .
\]
They satisfy the relation $d^{n+1}_{\Fh} \verk d^n_{\Fh} = 0$ so that we can form the leafwise cohomology of $\Fh$:
\[
H^n_{\Fh} (X) = \Ker d^n_{\Fh} / \Imm d^{n-1}_{\Fh} \; .
\]
For our purposes these invariants are actually too subtle. We therefore consider the reduced leafwise cohomology
\[
\oH^n_{\Fh} (X) = \Ker d^n_{\Fh} / \overline{\Imm d^{n-1}_{\Fh}} \; .
\]
Here the quotient is taken with respect to the topological closure of $\Imm d^{n-1}_{\Fh}$ in the natural Fr\'echet topology on $\Ah^n_{\Fh} (X)$. The reduced cohomologies are nuclear Fr\'echet spaces. Even if the leaves are dense, already $\oH^1_{\Fh} (X)$ can be infinite dimensional.

The cup product pairing induced by the exterior product of forms along the leaves turns $\oH^{\hullet}_{\Fh} (X)$ into a graded commutative $\oH^0_{\Fh} (X)$-algebra.

The Poincare Lemma extends to the foliation context and implies that
\[
H^n_{\Fh} (X) = H^n (X, \Rh) \; .
\]
Here $\Rh$ is the sheaf of smooth real valued functions which are locally constant on the leaves. In particular
\[
\oH^0_{\Fh} (X) = H^0_{\Fh} (X) = H^0 (X, \Rh) 
\]
consists only of constant functions if $\Fh$ contains a dense leaf.

For the torus foliation above with $\alpha \notin \Q$ we therefore have $\oH^0_{\Fh} (T^2) = \R$. Some Fourier analysis reveals that $\oH^1_{\Fh} (T^2) \cong \R$. The higher cohomologies vanish since almost by definition we have
\[
H^n_{\Fh} (X) = 0 \quad \mbox{for all} \; n > d = \dim \Fh \; .
\]
For a smooth map $f : X \to Y$ of foliated manifolds which maps leaves into leaves, continuous pullback maps
\[
f^* : \Ah^n_{\Fh_Y} (Y) \longrightarrow \Ah^n_{\Fh_X} (X)
\]
are defined for all $n$. They commute with $d_{\Fh}$ and respect the exterior product of forms. Hence they induce a continuous map of reduced cohomology algebras
\[
f^* : \oH^{\hullet}_{\Fh_Y} (Y) \longrightarrow \oH^{\hullet}_{\Fh_X} (X) \; .
\]
A (complete) flow is a smooth $\R$-action $\phi : \R \times X \to X , (t,x) \mapsto \phi^t (x)$. It is called $\Fh$-compatible if every diffeomorphism $\phi^t : X \to X$ maps leaves into leaves. If this is the case we obtain a linear $\R$-action $t \mapsto \phi^{t*}$ on $\oH^n_{\Fh} (X)$ for every $n$. Let
\[
\Theta : \oH^n_{\Fh} (X) \longrightarrow \oH^n_{\Fh} (X)
\]
denote the infinitesimal generator of $\phi^{t*}$:
\[
\Theta h = \lim_{t\to 0} \frac{1}{t} (\phi^{t*} h -h ) \; .
\]
The limit exists and $\Theta$ is continuous in the Fr\'echet topology. As $\phi^{t*}$ is an algebra endomorphism of the $\R$-algebra $\oH^{\hullet}_{\Fh} (X)$ it follows that $\Theta$ is an $\R$-linear derivation. Thus we have
\begin{equation}
  \label{eq:1}
  \Theta (h_1 \cup h_2) = \Theta h_1 \cup h_2 + h_1 \cup \Theta h_2
\end{equation}
for all $h_1 , h_2$ in $\oH^{\hullet}_{\Fh} (X)$.

For arbitrary foliations the reduced leafwise cohomology does not seem to have a good structure theory. For Riemannian foliations however the situation is much better. These foliations are characterized by the existence of a ``bundle-like'' metric $g$. This is a Riemannian metric whose geodesics are perpendicular to all leaves whenever they are perpendicular to one leaf. For example any one-codimensional foliation given by a closed one-form without singularities is Riemannian. 

The graded Fr\'echet space $\Ah^{\hullet}_{\Fh} (X)$ carries a canonical inner product:
\[
(\alpha , \beta) = \int_X \langle \alpha , \beta \rangle_{\Fh} \vol \; .
\]
Here $\langle , \rangle_{\Fh}$ is the Riemannian metric on $\Lambda^{\hullet} T^* \Fh$ induced by $g$ and $\vol$ is the volume form or density on $X$ coming from $g$. Let
\[
\Delta_{\Fh} = d_{\Fh} d^*_{\Fh} + d^*_{\Fh} d_{\Fh}
\]
denote the Laplacian using the formal adjoint of $d_{\Fh}$ on $X$. Since $\Fh$ is Riemannian the restriction of $\Delta_{\Fh}$ to any leaf $F$ is the Laplacian on $F$ with respect to the induced metric \cite{AK1} Lemma 3.2, i.e.
\[
(\Delta_{\Fh} \alpha) \, |_F = \Delta_F (\alpha \, |_F) \quad \mbox{for all} \; \alpha \in \Ah^{\hullet}_{\Fh} (X) \; .
\]
We now assume that $T \Fh$ is orientable. Via $g$ the choice of an orientation determines a volume form $\vol_{\Fh}$ in $\Ah^d_{\Fh} (X)$ and hence a Hodge $*$-operator
\[
*_{\Fh} : \Lambda^n T^*_x \Fh \silo \Lambda^{d-n} T^*_x \Fh \quad \mbox{for every} \; x \; \mbox{in} \; X \; .
\]
It is determined by the condition that
\[
v \wedge *_{\Fh} w = \langle v,w \rangle_{\Fh} \, \vol_{\Fh , x} \quad \mbox{for} \; v,w \; \mbox{in} \; \Lambda^{\hullet} T^*_x \Fh \; .
\]
These fibrewise star-operators induce the leafwise $*$-operator on forms:
\[
*_{\Fh} : \Ah^n_{\Fh} (X) \silo \Ah^{d-n}_{\Fh} (X) \; .
\]
We now list some important properties of leafwise cohomology.

{\bf Properties} Assume that $X$ is compact, $\Fh$ a $d$-dimensional oriented Riemannian foliation and $g$ a bundle-like metric for $\Fh$.

Then the natural map
\begin{equation}
  \label{eq:2}
  \Ker \Delta^n_{\Fh} \silo \oH^n_{\Fh} (X) \; ,\; \omega \longmapsto \omega \mod \overline{\Imm d^{n-1}_{\Fh}}
\end{equation}
is a topological isomorphism of Fr\'echet spaces. We denote its inverse by $\Hh$.

This result is due to \'Alvarez L\'opez and Kordyukov \cite{AK1}. It is quite deep since $\Delta_{\Fh}$ is only elliptic along the leaves so that the ordinary elliptic regularity theory does not suffice. For non-Riemannian foliations (\ref{eq:2}) does not hold in general \cite{DS1}. All the following results are consquences of this Hodge theorem.

The Hodge $*$-operator induces an isomorphism 
\[
*_{\Fh} : \Ker \Delta^n_{\Fh} \silo \Ker \Delta^{d-n}_{\Fh}
\]
since it commutes with $\Delta_{\Fh}$ up to sign. From (\ref{eq:2}) we therefore get isomorphisms for all $n$:
\begin{equation}
  \label{eq:3}
  *_{\Fh} : \oH^n_{\Fh} (X) \silo \oH^{d-n}_{\Fh} (X) \; .
\end{equation}
For the next property define the trace map
\[
\tr : \oH^d_{\Fh} (X) \longrightarrow \R
\]
by the formula
\[
\tr (h) = \int_X *_{\Fh} (h) \vol := \int_X *_{\Fh} (\Hh (h)) \vol \; .
\]
It is an isomorphism if $\Fh$ has a dense leaf. Note that for {\it any} representative $\alpha$ in the cohomology class $h$ we have
\[
\tr (h) = \int_X *_{\Fh} (\alpha) \vol \; .
\]
Namely $\alpha - \Hh (h) = d_{\Fh} \beta$ and
\begin{eqnarray*}
  \int_X *_{\Fh} (d_{\Fh} \beta) \vol & = & \pm \int_X d^*_{\Fh} (*_{\Fh} \beta) \vol \\
& = &\pm ( 1 , d^*_{\Fh} (*_{\Fh} \beta)) \\
& = & \pm ( d_{\Fh} (1) , *_{\Fh} \beta)\\
& = & 0 \; .
\end{eqnarray*}
Alternatively the trace functional is given by
\[
\tr (h) = \int_X \alpha \wedge *_{\perp} (1)
\]
where $*_{\perp} (1)$ is the transverse volume element for $g$ c.f. \cite{AK1}, \S\,3.

It is not difficult to see using (\ref{eq:2}) that we get a scalar product on $\oH^n_{\Fh} (X)$ for every $n$ by setting:
\begin{eqnarray}
  \label{eq:4}
  (h,h') & = & \tr (h \cup *_{\Fh} h') \\
& = & \int_X \langle \Hh (h) , \Hh (h') \rangle_{\Fh} \vol \; . \nonumber
\end{eqnarray}
It follows from this that the cup product pairing 
\begin{equation}
  \label{eq:5}
  \cup : \oH^n_{\Fh} (X) \times \oH^{d-n}_{\Fh} (X) \longrightarrow \oH^d_{\Fh} (X) \xrightarrow{tr} \R
\end{equation}
is non-degenerate.

Next we discuss the K\"unneth formula. Assume that $Y$ is another compact manifold with a Riemannian foliation $\Fh_Y$. Then the canonical map
\[
H^n_{\Fh_X} (X) \otimes H^m_{\Fh_Y} (Y) \longrightarrow H^{n+m}_{\Fh_X \times \Fh_Y} (X \times Y)
\]
induces a topological isomorphism \cite{M}:
\begin{equation}
  \label{eq:6}
  \oH^n_{\Fh_X} (X) \hat{\otimes} \oH^m_{\Fh_Y} (Y) \silo \oH^{n+m}_{\Fh_X \times \Fh_Y} (X \times Y) \; .
\end{equation}
Since the reduced cohomology groups are nuclear Fr\'echet spaces, it does not matter which topological tensor product is chosen in (\ref{eq:6}). The proof of this K\"unneth formula uses (\ref{eq:2}) and the spectral theory of the Laplacian $\Delta_{\Fh}$.

Before we deal with more specific topics let us mention that also Hodge--K\"ahler theory can be generalized. A complex structure on a foliation $\Fh$ is an almost complex structure $J$ on $T \Fh$ such that all restrictions $J \, |_F$ to the leaves are integrable. Then the leaves carry holomorphic structures which vary smoothly in the transverse direction. A foliation $\Fh$ with a complex structure $J$ is called K\"ahler if there is a hermitian metric $h$ on the complex bundle $T_c \Fh = (T \Fh , J)$ such that the K\"ahler form along the leaves
\[
\omega_{\Fh} = - \halb \imm h \in \Ah^2_{\Fh} (X)
\]
is closed. Note that for example any foliation by orientable surfaces can be given a K\"ahlerian structure by choosing a metric on $X$, c.f. \cite{MS} Lemma A.3.1. Let
\[
L_{\Fh} : \oH^n_{\Fh} (X) \longrightarrow \oH^{n+2}_{\Fh} (X) \; , \; L_{\Fh} (h) = h \cup [\omega_{\Fh}]
\]
denote the Lefschetz operator.

The following assertions are consequences of (\ref{eq:2}) combined with the classical Hodge--K\"ahler theory. See \cite{DS3} for details. Let $X$ be a compact orientable manifold and $\Fh$ a K\"ahlerian foliation with respect to the hermitian metric $h$ on $T_c \Fh$. Assume in addition that $\Fh$ is Riemannian. Then we have:
\begin{equation}
  \label{eq:7}
  \oH^n_{\Fh} (X) \otimes \C = \bigoplus_{p+q=n} H^{pq} \; , \quad \mbox{where} \; \overline{H^{pq}} = H^{qp} \; .
\end{equation}
Here $H^{pq}$ consists of those classes that can be represented by $(p,q)$-forms along the leaves. Moreover there are topological isomorphisms
\[
H^{pq} \cong \oH^q (X , \Omega^p_{\Fh})
\]
with the reduced cohomology of the sheaf of holomorphic $p$-forms along the leaves. 

Furthermore the Lefschetz operator induces isomorphisms
\begin{equation}
  \label{eq:8}
  L^i_{\Fh} : \oH^{d-i}_{\Fh} (X) \silo \oH^{d+i}_{\Fh} (X) \quad \mbox{for} \; 0 \le i \le d \; .
\end{equation}
Finally the space of primitive cohomology classes $\oH^n_{\Fh} (X)_{\prim}$ carries the structure of a polarizable $\ind \, \R$-Hodge structure of weight $n$.

After this review of important properties of the reduced leafwise cohomology of Riemannian foliations we turn to a specific result relating flows and cohomology.

\begin{theorem}
  \label{t21}
Let $X$ be a compact $3$-manifold and $\Fh$ a Riemannian foliation by surfaces with a dense leaf. Let $\phi^t$ be an $\Fh$-compatible flow on $X$ which is conformal on $T\Fh$ with respect to a metric $g$ on $T\Fh$ in the sense that for some constant $\alpha$ we have: 
\begin{equation}
  \label{eq:9}
  g (T_x \phi^t (v) , T_x \phi^t (w)) = e^{\alpha t} g (v,w) \; \mbox{for all} \; v,w \in T_x \Fh , x \in X \; \mbox{and} \; t \in \R \; .
\end{equation}
Then we have for the infinitesimal generator of $\phi^{t*}$ that:
\[
\Theta = 0 \; \mbox{on} \; \oH^0_{\Fh} (X) = \R \quad \mbox{and} \quad \Theta = \alpha \; \mbox{on} \; \oH^2_{\Fh} (X) \cong \R \; .
\]
On $\oH^1_{\Fh} (X)$ the operator $\Theta$ has the form
\[
\Theta = \frac{\alpha}{2} + S
\]
where $S$ is skew-symmetric with respect to the inner product $( , )$ above.
\end{theorem}

\begin{remarknn}
   For the bundle-like metric on $X$ required for the construction of $(,)$ we take any extension of the given metric on $T \Fh$ to a bundle-like metric on $TX$. Such extensions exist.
\end{remarknn}

\begin{proofof}
  {\bf \ref{t21}} Because we have a dense leaf, $\oH^0_{\Fh} (X) = H^0_{\Fh} (X)$ consists only of constant functions. On these $\phi^{t*}$ acts trivially so that $\Theta = 0$. Since
\[
*_{\Fh} : \R = \oH^0_{\Fh} (X) \silo \oH^2_{\Fh} (X)
\]
is an isomorphism and since
\[
\phi^{t*} (*_{\Fh} (1)) = e^{\alpha t} (*_{\Fh} 1)
\]
by conformality, we have $\Theta = \alpha$ on $\oH^2_{\Fh} (X)$.

For $h_1 , h_2$ in $\oH^1_{\Fh} (X)$ we find
\begin{equation}
  \label{eq:10}
  \alpha (h_1 \cup h_2) = \Theta (h_1 \cup h_2) \stackrel{\rm (\ref{eq:1})}{=} \Theta h_1 \cup h_2 + h_1 \cup \Theta h_2 \; .
\end{equation}
By conformality $\phi^{t*}$ commutes with $*_{\Fh}$ on $\oH^1_{\Fh} (X)$. Differentiating, it follows that $\Theta$ commutes with $*_{\Fh}$ as well. Since by definition we have
\[
(h, h') = \tr (h \cup *_{\Fh} h') \quad \mbox{for} \; h , h' \in \oH^1_{\Fh} (X) \; ,
\]
it follows from (\ref{eq:10}) that as desired:
\[
\alpha ( h , h' ) = ( \Theta h , h' ) + ( h , \Theta h' ) \; .
\]
\end{proofof}
\section{Dynamical Lefschetz trace formulas}
The formulas we want to consider in this section relate the compact orbits of a flow with the alternating sum of suitable traces on cohomology. A suggestive but non-rigorous argument of Guillemin \cite{Gu} later rediscovered by Patterson \cite{P} led to the following conjecture \cite{D2} \S\,3. Let $X$ be a compact manifold with a one-codimensional foliation $\Fh$ and an $\Fh$-compatible flow $\phi$. Assume that the fixed points and the periodic orbits of the flow are non-degenerate in the following sense: For any fixed point $x$ the tangent map $T_x \phi^t$ should have eigenvalues different from $1$ for all $t > 0$. For any closed orbit $\gamma$ of length $l (\gamma)$ and any $x \in \gamma$ and integer $k \neq 0$ the automorphism $T_x \phi^{kl (\gamma)}$ of $T_x X$ should have the eigenvalue $1$ with algebraic multiplicity one. Observe that the vector field $Y_{\phi}$ generated by the flow provides an eigenvector $Y_{\phi,x}$ for the eigenvalue $1$. 

Recall that the length $l (\gamma) > 0$ of $\gamma$ is defined by the isomorphism:
\[
\R / l (\gamma) \Z \silo \gamma \; , \; t \longmapsto \phi^t (x) \; .
\]
For a fixed point $x$ we set\footnote{This is different from the normalization in \cite{D2} \S\,3.}
\[
\varepsilon_x = \sgn \det (1 - T_x \phi^t \tei T_x \Fh) \; .
\]
This is independent of $t > 0$. For a closed orbit $\gamma$ and $k \in \Z \ohne 0$ set$^1$
\[
\varepsilon_{\gamma} (k) = \sgn \det (1 - T_x \phi^{kl (\gamma)} \tei T_x X / \R Y_{\phi,x}) = \sgn \det ( 1 - T_x \phi^{kl (\gamma)} \tei T_x \Fh) \; .
\]
It does not depend on the point $x \in \gamma$.

Finally let $\Dh' (J)$ denote the space of Schwartz distributions on an open subset $J$ of $\R$. 

\begin{conj}
  \label{t31}
For $X , \Fh$ and $\phi$ as above there exists a natural definition of a $\Dh' (\R^{>0})$-valued trace of $\phi^*$ on the reduced leafwise cohomology $\oH^{\hullet}_{\Fh} (X)$ such that in $\Dh' (\R^{>0})$ we have:
\begin{equation}
  \label{eq:11}
\sum\limits^{\dim \Fh}_{n=0} (-1)^n \Tr (\phi^* \tei \oH^n_{\Fh} (X)) = \sum\limits_{\gamma} l (\gamma) \sum\limits^{\infty}_{k=1} \varepsilon_{\gamma} (k) \delta_{kl (\gamma)} + \sum\limits_x \varepsilon_x |1 - e^{\kappa_x t}|^{-1} \; .
\end{equation}
\end{conj}

Here $\gamma$ runs over the closed orbits of $\phi$ which are not contained in a leaf and $x$ over the fixed points. For $a \in \R , \delta_a$ is the Dirac distribution in $a$ and $\kappa_x$ is defined by the action of $T_x \phi^t$ on the $1$-dimensional vector space $T_x X / T_x \Fh$. That action is multiplication by $e^{\kappa_x t}$ for some $\kappa_x \in \R$ and all $t$.

The conjecture is not known (except for $\dim X = 1$) if $\phi$ has fixed points. It may well have to be amended somewhat in that case. The analytic difficulty in the presence of fixed points lies in the fact that in this case $\Delta_{\Fh}$ has no chance to be transversally elliptic to the $\R$-action by the flow, so that the methods of transverse index theory do not apply directly. In the simpler case when the flow is everywhere transversal to $\Fh$, \'Alvarez L\'opez and Kordyukov have proved a  beautiful strengthening of the conjecture. Partial results were obtained by other methods in \cite{Laz}, \cite{DS2}. We now describe their result in a convenient way for our purposes:

\begin{punkt}
  \label{t32} \rm
Assume $X$ is a compact oriented manifold with a one codimensional foliation $\Fh$. Let $\phi$ be a flow on $X$ which is everywhere transversal to the leaves of $\Fh$. Then $\Fh$ inherits an orientation and it is Riemannian \cite{Go} III 4.4. Fixing a bundle-like metric $g$ the cohomologies $\oH^n_{\Fh} (X)$ acquire pre-Hilbert structures (\ref{eq:4}) and we can consider their Hilbert space completions $\hat{H}^n_{\Fh} (X)$. For every $t$ the linear operator $\phi^{t*}$ is bounded on $(\oH^n_{\Fh} (X) , \| \; \|)$ and hence can be continued uniquely to a bounded operator on $\hat{H}^n_{\Fh} (X)$ c.f. theorem \ref{t34}.
\end{punkt}

By transversality the flow has no fixed points. We assume that all periodic orbits are non-degenerate.

\begin{theorem}[\cite{AK2}]
  \label{t33}
Under the conditions of (\ref{t32}), for every test function $\varphi \in \Dh (\R) = C^{\infty}_0 (\R)$ the operator
\[
A_{\varphi} = \int_{\R} \varphi (t) \phi^{t*} \, dt
\]
on $\hat{H}^n_{\Fh} (X)$ is of trace class. Setting:
\[
\Tr (\phi^* \tei \oH^n_{\Fh} (X)) (\varphi) = \tr A_{\varphi}
\]
defines a distribution on $\R$. The following formula holds in $\Dh' (\R)$:
\begin{equation}
  \label{eq:12}
  \sum^{\dim \Fh}_{n=0} (-1)^n \Tr (\phi^* \tei \oH^n_{\Fh} (X)) = \chi_{\Co} (\Fh , \mu) \delta_0 + \sum_{\gamma} l (\gamma) \sum_{k \in \Z \ohne 0} \varepsilon_{\gamma} (k) \delta_{kl (\gamma)} \; .
\end{equation}
Here $\chi_{\Co} (\Fh , \mu)$ denotes Connes' Euler characteristic of the foliation with respect to the transverse measure $\mu$ corresponding to $*_{\perp} (1)$. (See \cite{MS}.)
\end{theorem}

It follows from the theorem that if the right hand side of (\ref{eq:12}) is non-zero, at least one of the cohomology groups $\oH^n_{\Fh} (X)$ must be infinite dimensional. Otherwise the alternating sum of traces would be a smooth function and hence have empty singular support.

By the Hodge isomorphism (\ref{eq:2}) one may replace cohomology by the spaces of leafwise harmonic forms. The left hand side of the dynamical Lefschetz trace formula then becomes the $\Dh' (\R)$-valued transverse index of the leafwise de Rham complex. Note that the latter is transversely elliptic for the $\R$-action $\phi^t$. Transverse index theory with respect to compact group actions was initiated in \cite{A}. A definition for non-compact groups of a transverse index was later given by H\"ormander \cite{Si} Appendix II. 

As far as we know the relation of (\ref{eq:12}) with transverse index theory in the sense of Connes--Moscovici still needs to be clarified.

Let us now make some remarks on the operators $\phi^{t*}$ on $\hat{H}^n_{\Fh} (X)$ in a more general setting:

\begin{theorem}
  \label{t34}
Let $\Fh$ be a Riemannian foliation on a compact manifold $X$ and $g$ a bundle like metric. As above $\hat{H}^n_{\Fh} (X)$ denotes the Hilbert space completion of $\oH^n_{\Fh} (X)$ with respect to the scalar product (\ref{eq:4}). Let $\phi^t$ be an $\Fh$-compatible flow. Then the linear operators $\phi^{t*}$ on $\oH^n_{\Fh} (X)$ induce a strongly continuous operator group on $\hat{H}^n_{\Fh} (X)$. In particular the infinitesimal generator $\Theta$ exists as a closed densely defined operator. On $\oH^n_{\Fh} (X)$ it agrees with the infinitesimal generator in the Fr\'echet topology defined earlier. There exists $\omega > 0$ such that the spectrum of $\Theta$ lies in $- \omega \le \RRe s \le \omega$. If the operators $\phi^{t*}$ are orthogonal then $T = - i \Theta$ is a self-adjoint operator on $\hat{H}^n_{\Fh} (X) \otimes \C$ and we have
\[
\phi^{t*} = \exp t \Theta = \exp it T
\]
in the sense of the functional calculus for (unbounded) self-adjoint operators on Hilbert spaces. 
\end{theorem}

{\bf Sketch of proof} Estimates show that $\| \phi^{t*} \|$ is locally uniformly bounded in $t$ on $\oH^n_{\Fh} (X)$. Approximating $h \in \hat{H}^n_{\Fh} (X)$ by $h_{\nu} \in \oH^n_{\Fh} (X)$ one now shows as in the proof of the Riemann--Lebesgue lemma that the function $t \mapsto \phi^{t*} h$ is continuous at zero, hence everywhere. Thus $\phi^{t*}$ defines a strongly continuous group on $\hat{H}^n_{\Fh} (X)$. The remaining assertions follow from semigroup theory \cite{DSch}, Ch. VIII, XII, and in particular from the theorem of Stone.

We now combine theorems \ref{t21}, \ref{t33} and \ref{t34} to obtain the following corollary:

\begin{cor}
  \label{t35}
Let $X$ be a compact $3$-manifold with a foliation $\Fh$ by surfaces having a dense leaf. Let $\phi^t$ be a non-degenerate $\Fh$-compatible flow which is everywhere transversal to $\Fh$. Assume that $\phi^t$ is conformal (\ref{eq:9}) with respect to a metric $g$ on $T\Fh$. Then $\Theta$ has pure point spectrum $\Sp^1 (\Theta)$ on $\hat{H}^1_{\Fh} (X)$ which is discrete in $\R$ and we have the following equality of distributions on $\R$:
\begin{equation}
  \label{eq:13}
  1 - \sum_{\rho \in \Sp^1 (\Theta)} e^{t\rho} + e^{t\alpha} = \chi_{\Co} (\Fh , \mu) \delta_0 + \sum_{\gamma} l (\gamma) \sum_{k \in \Z \ohne 0} \varepsilon_{\gamma} (k) \delta_{k l (\gamma)} \; .
\end{equation}
In the sum the $\rho$'s appear with their geometric multiplicities. All $\rho \in \Sp^1 (\Theta)$ have $\RRe \rho = \frac{\alpha}{2}$.
\end{cor}

{\bf Remarks} 1) Here $e^{t\rho} , e^{t \alpha}$ are viewed as distributions so that evaluated on a test function $\varphi \in \Dh (\R)$ the formula reads:
\begin{equation}
  \label{eq:14}
\Phi (0) - \sum\limits_{\rho \in \Sp^1 (\Theta)} \Phi (\rho) + \Phi (\alpha) = \chi_{\Co} (\Fh , \mu) \varphi (0) + \sum\limits_{\gamma} l (\gamma) \sum\limits_{k \in \Z \ohne 0} \varepsilon_{\gamma} (k) \varphi (k l (\gamma)) \; .
\end{equation}
Here we have put
\[
\Phi (s) = \int_{\R} e^{ts} \varphi (t) \, dt \; .
\]
2) Actually the conditions of the corollary force $\alpha = 0$ i.e. the flow must be isometric with respect to $g$. We have chosen to leave the $\alpha$ in the fomulation since there are good reasons to expect the corollary to generalize to more general phase spaces $X$ than manifolds, where $\alpha \neq 0$ becomes possible i.e. to Sullivan's generalized solenoids. More on this in section 5.\\
3) One can show that the group generated by the lengths of closed orbits is a finitely generated subgroup of $\R$ under the assumptions of the corollary. In order to achieve an infinitely generated group the flow must have fixed points.

\begin{proofof}
  {\bf \ref{t35}}
By \ref{t21}, \ref{t33} we need only show the equation
\begin{equation}
  \label{eq:15}
  \Tr (\phi^{t*} \tei \oH^1_{\Fh} (X)) = \sum_{\rho \in \Sp^1 (\Theta)} e^{t \rho}
\end{equation}
and the assertions about the spectrum of $\Theta$. As in the proof of \ref{t21} one sees that on $\hat{H}^1_{\Fh} (X)$ we have
\[
( \phi^{t*} h , \phi^{t*} h') = e^{\alpha t} (h , h') \; .
\]
Hence $e^{-\frac{\alpha}{2} t} \phi^{t*}$ is orthogonal and by the theorem of Stone
\[
T = -iS
\]
is selfadjoint on $\hat{H}^1_{\Fh} (X) \otimes \C$, if $\Theta = \frac{\alpha}{2} + S$. Moreover
\[
e^{-\frac{\alpha}{2} t} \phi^{t*} = \exp it T \; ,
\]
so that
\begin{equation}
  \label{eq:16}
  \phi^{t*} = \exp t \Theta \; .
\end{equation}
In \cite{DS2} proof of 2.6, for isometric flows the relation
\[
-\Theta^2 = \Delta^1 \, |_{\ker \Delta^1_{\Fh}}
\]
was shown. Using the spectral theory of the ordinary Laplacian $\Delta^1$ on $1$-forms it follows that $\Theta$ has pure point spectrum with finite multiplicities on $\oH^1_{\Fh} (X) \cong \ker \Delta^1_{\Fh}$ and that $\Sp^1 (\Theta)$ is discrete in $\R$. Alternatively, without knowing $\alpha = 0$, that proof gives:
\[
- \left( \Theta - \frac{\alpha}{2} \right)^2 = \Delta^1 \, |_{\ker \Delta^1_{\Fh}} \; .
\]
This also implies the assertion on the spectrum of $\Theta$ on $\hat{H}^1_{\Fh} (X)$. Now (\ref{eq:15}) follows from (\ref{eq:16}) and the fact that the operators $A_{\varphi}$ are of trace class.
\end{proofof}

\begin{remarknn}
In more general situations where $\Theta$ may not have a pure point spectrum of $\hat{H}^1_{\Fh} (X)$ but where $e^{-\frac{\alpha}{2}} \phi^{t*}$ is still orthogonal, we obtain:
\[
\langle \Tr (\phi^* \tei \oH^1_{\Fh} (X)) , \varphi \rangle = \sum_{\rho \in \Sp^1 (\Theta)_{\mathrm{point}}} \Phi (\rho) + \int^{\frac{\alpha}{2} + i \infty}_{\frac{\alpha}{2} - i \infty} \Phi (\lambda) m (\lambda) \, d \lambda
\]
where $m (\lambda) \ge 0$ is the spectral density function of the continuous part of the spectrum of $\Theta$. 
\end{remarknn}
\section{Comparison with the ``explicit formulas'' in analytic number theory}
Consider a number field $K / \Q$. The explicit formulas in analytic number theory relate the primes of $K$ to the non-trivial zeroes of the Dedekind zeta function $\zeta_K (s)$ of $K$.

\begin{theorem}
  \label{t41}
For $\varphi \in \Dh (\R)$ define $\Phi (s)$ as in the preceeding section. Then the following fomula holds:
\begin{eqnarray}
  \label{eq:17}
  \lefteqn{\Phi (0) - \sum_{\rho} \Phi (\rho) + \Phi (1) = - \log |d_{K / \Q}| \varphi (0)} \nonumber \\
& & + \sum_{\ep \nmid \infty} \log N \ep \left( \sum_{k \ge 1} \varphi (k \log N \ep) + \sum_{k \le -1} N \ep^k \varphi (k \log N \ep) \right) \nonumber \\
& & + \sum_{\ep \tei \infty} W_{\ep} (\varphi) \; .
\end{eqnarray}
\end{theorem}

Here $\rho$ runs over the non-trivial zeroes of $\zeta_K (s)$ i.e. those that are contained in the critical strip $0 < \RRe s < 1$. Moreover $\ep$ runs over the places of $K$ and $d_{K / \Q}$ is the discriminant of $K$ over $\Q$. For $\ep \tei \infty$ the $W_{\ep}$ are distributions which are determined by the $\Gamma$-factor at $\ep$. If $\varphi$ has support in $\R^{> 0}$ then
\[
W_{\ep} (\varphi) = \int^{\infty}_{-\infty} \frac{\varphi (t)}{1 - e^{\kappa_{\ep} t}} \, dt
\]
where $\kappa_{\ep} = -1$ if $\ep$ is complex and $\kappa_{\ep} = -2$ if $\ep$ is real. If $\varphi$ has support on $\R^{< 0}$ then
\[
W_{\ep} (\varphi) = \int^{\infty}_{-\infty} \frac{\varphi (t)}{1 - e^{\kappa_{\ep} |t|}} \, e^t \, dt \; .
\]
There are different ways to write $W_{\ep}$ on all of $\R$ but we will not discuss this here. See for example \cite{Ba} which also contains a proof of the theorem for much more general test functions.

Formula (\ref{eq:17}) implies the following equality of distributions on $\R^{> 0}$:
\begin{equation}
  \label{eq:18}
  1 - \sum_{\rho} e^{t\rho} + e^t = \sum_{\ep \nmid \infty} \log N\ep \sum^{\infty}_{k=1} \delta_{k \log N \ep} + \sum_{\ep \tei \infty} (1 - e^{\kappa_{\ep} t})^{-1} \; .
\end{equation}
This fits rather nicely with formula (\ref{eq:11}) and suggests the following analogies:
\begin{center}
  \begin{tabular}{llp{8cm}}
$\spec \eo_K \cup \{ \ep \tei \infty \}$ & \quad $\ent$ \quad & $3$-dimensional dynamical system $(X , \phi^t)$ with a one-codimensional foliation $\Fh$ satisfying the conditions of conjecture \ref{t31} \\
finite place $\ep$ & \quad $\ent$ \quad & closed orbit $\gamma = \gamma_{\ep}$ not contained in a leaf and hence transversal to $\Fh$ such that $l (\gamma_{\ep}) = \log N \ep$ and $\varepsilon_{\gamma_{\ep}} (k) = 1$ for all $k \ge 1$.\\
infinite place $\ep$ & \quad $\ent$ \quad & fixed point $x_{\ep}$ such that $\kappa_{x_{\ep}} = \kappa_{\ep}$ and \newline
$\varepsilon_{x_{\ep}} = 1$.
  \end{tabular}
\end{center}

In order to understand number theory more deeply in geometric terms it would be very desirable to find a system $(X , \phi^t , \Fh)$ which actually realizes this correspondence. For this the class of compact $3$-manifolds as phase spaces has to be generalized as will become clear from the following discussion.

Formula (\ref{eq:17}) can be written equivalently as an equality of distributions on $\R$:
\begin{eqnarray}
  \label{eq:19}
  \lefteqn{1 - \sum_{\rho} e^{t\rho} + e^t  =  - \log |d_{K / \Q}| \delta_0} \nonumber \\
 & &  + \sum_{\ep \nmid \infty} \log N\ep\left( \sum_{k \ge 1} \delta_{k \log N \ep} + \sum_{k \le 1} N \ep^{k} \delta_{k \log N \ep} \right) \nonumber \\
& & + \sum_{\ep \tei \infty} W_{\ep} \; .
\end{eqnarray}

Let us compare this with formula (\ref{eq:13}) in Corollary \ref{t35}. This corollary is the best result yet on the dynamical side but still only a first step since it does not allow for fixed points which as we have seen must be expected for dynamical systems of relevance for number fields.

Ignoring the contributions $W_{\ep}$ from the infinite places for the moment we are suggested that
\begin{equation}
  \label{eq:20}
  - \log |d_{K/ \Q}| \ent \chi_{Co} (\Fh , \mu) \; .
\end{equation}
There are two nice points about this analogy. Firstly there is the following well known fact due to Connes: 

\begin{fact}
  \label{t42}
Let $\Fh$ be a foliation of a compact $3$-manifold by surfaces such that the union of the compact leaves has $\mu$-measure zero, then
\[
\chi_{Co} (\Fh , \mu) \le 0 \; .
\]
\end{fact}
Namely the non-compact leaves are known to be complete in the induced metric. Hence they carry no non-zero harmonic $L^2$-functions, so that Connes' $0$-th Betti number $\beta_0 (\Fh , \mu) = 0$. Since $\beta_2 (\Fh , \mu) = \beta_0 (\Fh , \mu)$ it follows that
\begin{eqnarray*}
  \chi_{Co} (\Fh , \mu) & = & \beta_0 (\Fh , \mu) - \beta_1 (\Fh , \mu) + \beta_2 (\Fh , \mu) \\
& = & - \beta_1 (\Fh , \mu) \le 0 \; .
\end{eqnarray*}
The reader will have noticed that in accordance with \ref{t42} the left hand side of (\ref{eq:20}) is negative as well:
\[
- \log |d_{K / \Q}| \le 0 \quad \mbox{for all} \;  K / \Q \; .
\]
The second nice point about (\ref{eq:20}) is this. The bundle-like metric $g$ which we have chosen for the definition of $\Delta_{\Fh}$ and of $\chi_{Co} (\Fh , \mu)$ induces a holomorphic structure of $\Fh$ \cite{MS}, Lemma A\,3.1. The space $X$ is therefore foliated by Riemann surfaces. Let $\chi_{Co} (\Fh , \Oh , \mu)$ denote the holomorphic Connes Euler characteristic of $\Fh$ defined using $\Delta_{\overline{\partial}}$-harmonic forms on the leaves instead of $\Delta$-harmonic ones. According to Connes' Riemann--Roch Theorem \cite{MS} Cor. A.\,2.3, Lemma A\,3.3 we have:
\[
\chi_{Co} (\Fh , \Oh , \mu) = \halb \chi_{Co} (\Fh , \mu) \; .
\]
Therefore $\chi_{Co} (\Fh , \Oh , \mu)$ corresponds to $-\log \sqrt{|d_{K / \Q}|}$. 

For completely different reasons this number is defined in Arakelov theory as the Arakelov Euler characteristic of $\overline{\spec \eo_K} = \spec \eo_K \cup \{ \ep \tei \infty \}$:
\begin{equation}
  \label{eq:21}
  \chi_{Ar} (\Oh_{\overline{\spec \eo_K}}) = - \log {\textstyle \sqrt{|d_{K / \Q}|}} \; .
\end{equation}
See \cite{N} for example. Thus we see that
\[
\chi_{Ar} (\Oh_{\overline{\spec \eo_K}}) \quad \mbox{corresponds to} \; \chi_{Co} (\Fh , \Oh , \mu) \; .
\]
It would be very desirable of course to understand Arakelov Euler characteristics in higher dimensions even conjecturally in terms of Connes' holomorphic Euler characteristics. Note however that Connes' Riemann--Roch theorem in higher dimensions does not involve the $R$-genus appearing in the Arakelov Riemann--Roch theorem. The ideas of Bismut \cite{Bi} may be relevant in this connection. He interpretes the $R$-genus in a natural way via the geometry of loop spaces. 

Further comparison of formulas (\ref{eq:13}) and (\ref{eq:19}) shows that in a dynamical system corresponding to number theory we must have $\alpha = 1$. This means that the flow $\phi^{t*}$ would act by multiplication  with $e^t$ on the one-dimensional space $\oH^2_{\Fh} (X)$. As explained before this would be the case if $\phi^t$ were conformal on $T\Fh$ with factor $e^t$:
\begin{equation}
  \label{eq:22}
  g (T_x \phi^t (v) , T_x \phi^t (w)) = e^t g (v,w) \quad \mbox{for all} \; v,w \in T_x \Fh \; .
\end{equation}
However as mentioned before, this is not possible in the manifold setting of corollary \ref{t35} which actually implies $\alpha = 0$. \\
An equally important difference between formulas (\ref{eq:13}) and (\ref{eq:19}) is between the coefficients of $\delta_{kl (\gamma)}$ and of $\delta_{k \log N\ep}$ for $k \le -1$. In the first case it is $\pm 1$ whereas in the second it is $N\ep^k = e^{k \log N \ep}$ which corresponds to $e^{kl (\gamma)}$. 

Thus it becomes vital to find phase spaces $X$ more general than manifolds for which the analogue of corollary \ref{t35} holds and where $\alpha \neq 0$ and in particular $\alpha = 1$ becomes possible. In the new context the term $\varepsilon_{\gamma} (k) \delta_{k l (\gamma)}$ for $k \le -1$ in formula (\ref{eq:13}) should become $\varepsilon_{\gamma} (k) e^{\alpha k l (\gamma)} \delta_{kl (\gamma)}$. The next section is devoted to a discussion of certain laminated spaces which we propose as possible candidates for this goal.
\section{Remarks on dynamical Lefschetz trace formulas on laminated spaces}
In this section we extend the previous discussion to more general phase spaces than manifolds. The class of spaces we have in mind are the foliated spaces with totally disconnected transversals in the sense of \cite{MS}. We will call them laminated spaces for short. They also go by the name of (generalized) solenoids c.f. \cite{Su}.

\begin{punkt}
  \label{t51} \rm
An $a$-dimensional laminated space is a second countable metrizable topological space $X$ which is locally homeomorphic to the product of a non-empty open subset of $\R^a$ with a totally disconnected space. Then $a$ is the topological dimension of $X$.

Transition functions between local charts $\varphi_1$ and $\varphi_2$ have the following form locally:
\begin{equation}
  \label{eq:23}
  \varphi_2 \verk \varphi^{-1}_1 (x,y) = (F_1 (x,y) , F_2 (y)) \; .
\end{equation}
Here $x,y$ denote the euklidean resp. totally disconnected components. This is due to the fact, that continuous functions from connected subsets of $\R^a$ into a totally disconnected space are constant.

Because of (\ref{eq:23}) the inverse images $\varphi^{-1} (\, , *)$ patch together, to give a partition $\Lh$ of $X$ into $a$-dimensional topological manifolds. These {\it leaves} of the laminated space $X$ are exactly the path components of $X$. The classical solenoid
\begin{equation}
  \label{eq:24}
  \Sa^1_p = \R \times_{\Z} \Z_p = \lim_{\leftarrow} (\ldots \to \R / \Z \xrightarrow{p} \R / \Z \xrightarrow{p} \ldots )
\end{equation}
is an example of a compact connected one-dimensional laminated space with dense leaves homeomorphic to the real line.

A $C^{\infty , 0}$-structure on a laminated space is a maximal atlas of local charts whose transition functions are smooth in the euklidean component and continuous in the totally disconnected one. Furthermore all derivatives in the euklidean directions should be continuous in all components. The leaves then become $a$-dimensional smooth manifolds.

$C^{\infty , 0}$-laminated spaces are examples of foliated spaces in the sense of \cite{MS} Def. 2.1.

A stronger structure that may exist on a laminated space was introduced by Sullivan \cite{}. A $C^{\infty , \infty}$- or $TLC$-structure on a laminated space $X$ is given by a maximal atlas whose transition functions are smooth in the euklidean component and uniformly locally constant in the totally disconnected one. That is, locally they have the form:
\begin{equation}
  \label{eq:25}
  \varphi_2 \verk \varphi^{-1}_1 (x,y) = (F_1 (x) , F_2 (y))
\end{equation}
with $F_1$ smooth and $F_2$ locally constant. Every $C^{\infty , \infty}$-structure gives rise to a $C^{\infty , 0}$-structure. It is clear that $\Sa^1_p$ is naturally a $C^{\infty , \infty}$-laminated space.

For a $C^{\infty , 0}$-laminated space $X$ let $TX = T \Lh$ denote its tangent bundle in the sense of \cite{MS} p. 43. For a point $x \in X$, the fibre $T_x X$ is the ordinary tangent space to the leaf through $x$. A Riemannian metric on $X$ is one on $TX$. Morphisms between $C^{\infty , 0}$-laminated spaces are continuous maps which induce smooth maps between the leaves of the lamination. They induce morphisms of tangent bundles.

The two most prominent places in mathematics where laminated spaces occur naturally are in number theory e.g. as adelic points of algebraic groups and in the theory of dynamical systems as attractors.
\end{punkt}

\begin{punkt}
  \label{t52} \rm
We now introduce foliations of laminated spaces. Let $X$ be an $a$-dimensional laminated space. For our purposes, a foliation $\Fh$ of $X$ by laminated spaces is a partition of $X$ into $d$-dimensional laminated spaces. The foliation is supposed to be locally trivial with euclidean transversals. More precisely $\Fh$ is given by a maximal atlas of local charts on $X$
\[
\varphi : U \silo V_1 \times V_2 \times Y
\]
with $V_1 \subset \R^d , V_2 \subset \R^{a-d}$ open and $Y$ totally disconnected, having the following property:

The transition maps have the form:
\[
\varphi_2 \verk \varphi^{-1}_1 (x_1 , x_2 , y) = (G_1 (x_1 , x_2 , y) , G_2 (x_2 , y) , G_3 (y))
\]
where $G_1 , G_2$ are smooth in the $x_1 , x_2$ components and $G_3$ and all $\partial^{\alpha_1 , \alpha_2}_{x_1 , x_2} G_1$ and $\partial^{\alpha_2}_{x_2} G_2$ are continuous. Setting $x = (x_1 , x_2)$, 
\[
F_1 (x,y) = (G_1 (x_1 , x_2 , y_2) , G_2 (x_2 , y))
\]
and $F_2 = G_2$ this induces a $C^{\infty , 0}$-structure on $X$ which is supposed to agree with the given one.

The leaves of $\Fh$ are obtained by patching together the sets $\varphi^{-1} (V_1 \times \{ * \} \times Y )$. They are $d$-dimensional $C^{\infty , 0}$-laminated spaces with local transition functions given by:
\[
(x_1 , y) \longmapsto (G_1 (x_1 , * , y) , G_3 (y)) \; .
\]
Their leaves are $d$-dimensional manifolds which foliate the $a$-dimensional manifolds which occur as leaves of the lamination on $X$. Thus $X$ is partitioned into $a$-dimensional manifolds each of which carries a $d$-dimensional foliation in the usual sense.

Besides $\Lh$ and $\Fh$ there is a third foliated structure denoted $\Fh \Lh$ on $X$. The space $X$ is foliated with leaves the $d$-dimensional manifolds that occur as path components of the $\Fh$-leaves. Here the transverse space is of the form
\begin{center}
  open subspace of $\R^{a-d} \times$ totally disconnected.
\end{center}
The local transition maps are given by:
\[
(x_1 , z) \longmapsto (H_1 (x_1 , z) , H_2 (z))
\]
with $z = (x_2 , y) , H_1 (x_1 , z) = G_1 (x_1 , x_2 , y)$ and $H_2 (z) = (G_2 (x_2 , y) , G_3 (y))$. 

Of the three foliated structures $\Lh , \Fh$ and $\Fh \Lh$ on $X$ the first and the last fit into the context of \cite{MS} but not the second. 
\end{punkt}

\begin{punkt}
  \label{t53} \rm
We now turn to cohomology. The foliation $\Fh\Lh$ gives rise to the integrable rank $d$ subbundle $T\Fh\Lh$ of $TX = T \Lh$. These bundles are tangentially smooth in the sense of \cite{MS}, p. 43 with respect to $\Lh$. The leafwise cohomology of $X$ along $\Fh \Lh$ is by definition the cohomology of the sheaf $\Rh_{\Fh \Lh}$ of real valued smooth functions on $X$ wich are locally constant along the $\Fh \Lh$ leaves. Here a continuous function or form on an open subset of $X$ is called smooth if its restrictions to the $\Lh$-leaves are smooth. The sheaf $\Rh_{\Fh \Lh}$ is resolved by the de Rham complex of smooth forms along $\Fh \Lh$. Using the natural Fr\'echet topology on the spaces of global differential forms, one defines the reduced version $\oH^{\hullet}_{\Fh\Lh} (X)$ of leafwise cohomology $H^{\hullet}_{\Fh\Lh} (X)$ as its maximal Hausdorff quotient. By $H^{\hullet}_{\Fh} (X)$ we denote the cohomology of the sheaf $\Rh_{\Fh}$ of real valued smooth functions on $X$ which are locally constant along the $\Fh$-leaves.
\end{punkt}

\begin{punkt}
  \label{t54} \rm
A flow $\phi$ on $X$ is a continuous $\R$-action such that the induced $\R$-actions on the leaves of $\Lh$ are smooth. It respects $\Fh$ if every $\phi^t$ maps leaves of $\Fh$ into leaves of $\Fh$. It follows that $\phi^t$ maps $\Fh \Lh$-leaves into $\Fh\Lh$-leaves. Thus $(X , \Fh , \phi^t)$ is partitioned into the foliated dynamical systems $(L, \Fh_L , \phi^t \, |_L)$ for $L \in \Lh$. Here
\[
\Fh_L = \Fh \Lh \, |_L = \{ S \in \Fh \Lh \tei S \subset L \} \; .
\]
Any $\Fh$-compatible flow $\phi^t$ induces pullback actions $\phi^{t*}$ on $H^{\hullet}_{\Fh\Lh} (X)$ and $\oH^{\hullet}_{\Fh\Lh} (X)$.
\end{punkt}

\begin{punkt}
  \label{t55} \rm
We now state as a working hypotheses a generalization of the conjectured dynamical trace formula \ref{t31}. We allow the phase space to be a laminated space. Moreover we extend the formula to an equality of distributions on $\R^*$ instead of $\R^{>0}$. After checking various compatibilities we state a case where our working hypotheses can be proved and give a number theoretical example.
\end{punkt}

\begin{punkt} {\bf Working hypotheses:}
  \label{t56} \rm
Let $X$ be a compact $C^{\infty , 0}$-laminated space with a one-codimensional foliation $\Fh$ and an $\Fh$-compatible flow $\phi$. Assume that the fixed points and the periodic orbits of the flow are non-degenerate. Then there exists a natural definition of a $\Dh' (\R^*)$-valued trace of $\phi^{t*}$ on $\oH^{\hullet}_{\Fh\Lh} (X)$ such that in $\Dh' (\R^*)$ we have:
\begin{eqnarray}
  \label{eq:26}
\hspace*{0.5cm}  \lefteqn{\sum^{\dim \Fh}_{n=0} (-1)^n \Tr (\phi^* \tei \oH^n_{\Fh\Lh} (X)) = }\\
& & \sum_{\gamma} l (\gamma) \left( \sum_{k \ge 1} \varepsilon_{\gamma} (k) \delta_{k l (\gamma)} + \sum_{k \le -1} \varepsilon_{\gamma} (|k|) \det (-T_x \phi^{kl (\gamma)} \tei T_x \Fh) \delta_{kl (\gamma)} \right) \nonumber \\
& & + \sum_x W_x \; . \nonumber
\end{eqnarray}
Here $\gamma$ runs over the closed orbits not contained in a leaf and in the sums over $k$'s any point $x \in \gamma$ can be chosen. The second sum runs over the fixed points $x$ of the flow. The distributions $W_x$ on $\R^*$ are given by:
\[
W_x \, |_{\R^{> 0}} = \varepsilon_x \, |1 - e^{\kappa_x t}|^{-1}
\]
and
\[
W_x \, |_{\R^{< 0}} = \varepsilon_x \det (-T_x \phi^t \tei T_x \Fh) \, |1 - e^{\kappa_x |t|}|^{-1} \; .
\]
\end{punkt}

\begin{remarks} \label{t57}
  \rm {\bf 0)} It may actually be better to use a version of foliation cohomology where transversally forms are only supposed to be locally $L^2$ instead of being continuous.\\
{\bf 1)} In the situation described in \ref{t58} below the working hypotheses can be proved if $\oH^n_{\Fh\Lh} (X)$ is replaced by $H^n_{\Fh} (X)$, Theorem \ref{t59}. In those cases there are no fixed points, only closed orbits. Thus Theorem \ref{t59} dictated only the coefficients of $\delta_{kl (\gamma)}$ for $k \in \Z \ohne 0$, but not the contributions $W_x$ from the fixed points. \\
{\bf 2)} The coefficients of $\delta_{kl (\gamma)}$ for $k \in \Z \ohne 0$ can be written in a uniform way as follows. They are equal to:
\begin{equation}
  \label{eq:27}
  \frac{\det (1 - T_x \phi^{kl (\gamma)} \tei T_x \Fh)}{|\det (1 - T_x \phi^{|k| l (\gamma)} \tei T_x X / \R  Y_{\phi , x})|} = \frac{\det (1 - T_x \phi^{kl (\gamma)} \tei T_x \Fh)}{| \det (1 - T_x \phi^{|k| l (\gamma)} \tei T_x \Fh)|} \; .
\end{equation}
Here $x$ is any point on $\gamma$. Namely, for $k \ge 1$ this equals $\varepsilon_{\gamma} (k)$ whereas for $k \le -1$ we obtain
\begin{equation}
  \label{eq:28}
  \varepsilon_{\gamma} (|k|) \det (- T_x \phi^{kl (\gamma)} \tei T_x \Fh) = \varepsilon_{\gamma} (k) \, |\det (T_x \phi^{kl (\gamma)} \tei T_x \Fh)| \; .
\end{equation}
The expression on the left hand side of (\ref{eq:27}) motivated our conjecture about the contributions on $\R^*$ from the fixed points $x$. Since $Y_{\phi , x} = 0$, they should be given by:
\[
\frac{\det (1 - T_x \phi^t \tei T_x \Fh)}{| \det (1 - T_x \phi^{|t|} \tei T_x X)|} \overset{!}{=} W_x \; .
\]
{\bf 3)} One can prove that in the manifold setting of theorem \ref{t33} we have
\[
|\det (T_x \phi^{kl (\gamma)} \tei T_x \Fh)| = 1 \; .
\]
By (\ref{eq:26}), our working hypotheses \ref{t56} is therefore compatible with formula (\ref{eq:12}). Compatibility with conjecture \ref{t31} is clear.\\
{\bf 4)} We will see below that in our new context metrics $g$ on $T\Fh$ can exist for which the flow has the conformal behaviour (\ref{eq:22}). Assuming we are in such a situation and that $\Fh$ is $2$-dimensional, we have:
\[
|\det (T_x \phi^{kl (\gamma)} \tei T_x \Fh)| = e^{kl (\gamma)} \quad \mbox{for} \; x \in \gamma , k \in \Z
\]
and
\[
|\det (T_x \phi^t \tei T_x \Fh)| = e^t \quad \mbox{for a fixed point} \; x \; .
\]
In the latter case, we even have by continuity:
\[
\det (T_x \phi^t \tei T_x \Fh) = e^t \; ,
\]
the determinant being positive for $t = 0$. Hence by (\ref{eq:28}) the conjectured formula (\ref{eq:26}) reads as follows in this case:
\begin{eqnarray}
  \label{eq:29}
  \lefteqn{\sum^2_{n=0} (-1)^n \Tr (\phi^* \tei \oH^n_{\Fh\Lh} (X))} \\
& = & \sum_{\gamma} l (\gamma) \left( \sum_{k \ge 1} \varepsilon_{\gamma} (k) \delta_{kl (\gamma)} + \sum_{k \le -1} \varepsilon_{\gamma} (k) e^{kl (\gamma)} \delta_{kl (\gamma)} \right) \nonumber \\
& & + \sum_x W_x \; . \nonumber
\end{eqnarray}
Here:
\[
W_x \, |_{\R^{> 0}} = \varepsilon_x \, |1 - e^{\kappa_x t}|^{-1}
\]
and
\[
W_x \, |_{\R^{< 0}} = \varepsilon_x e^t \, |1 - e^{\kappa_x |t|}|^{-1} \; .
\]
This fits perfectly with the explicit formula (\ref{eq:19}) if all $\varepsilon_{\gamma_{\ep}} (k) = 1$ and $\varepsilon_{x_{\ep}} = 1$. Namely if $l (\gamma_{\ep}) = \log N \ep$ for $\ep \nmid \infty$ and $\kappa_{x_{\ep}} = \kappa_{\ep}$ for $\ep \tei \infty$, then we have:
\[
e^{kl (\gamma_{\ep})} = e^{k \log N\ep} = N \ep^k \quad \mbox{for finite places} \; \ep
\]
and
\[
W_{x_{\ep}} = W_{\ep} \quad \mbox{on} \; \R^* \; \mbox{for the infinite places} \; \ep \; .
\]
{\bf 5)} In the setting of the preceeding remark the automorphisms
\[
e^{-\frac{k}{2} l (\gamma)} T_x \phi^{kl (\gamma)} \quad \mbox{of} \; T_x \Fh \quad \mbox{for} \; x \in \gamma
\]
respectively
\[
e^{-\frac{t}{2}} T_x \phi^t \quad \mbox{of} \; T_x \Fh \quad \mbox{for a fixed point} \; x
\]
are orthogonal automorphisms. For a real $2 \times 2$ orthogonal determinant $O$ with $\det O = -1$ we have:
\[
\det (1 - uO) = 1 - u^2 \; .
\]
The condition $\varepsilon_{\gamma} (k) = +1$ therefore implies that $\det (T_x \phi^{kl (\gamma)} \tei T_x \Fh)$ is positive for $k \ge 1$ and hence for all $k \in \Z$. The converse is also true. For a fixed point we have already seen directly that $\det (T_x \phi^t \tei T_x \Fh)$ is positive for all $t \in \R$. Hence we have the following information. 

{\bf Fact} In the situation of the preceeding remark, $\varepsilon_k (\gamma) = +1$ for all $k \in \Z \ohne 0$ if and only if on $T_x \Fh$ we have:
\[
T_x \phi^{kl (\gamma)} = e^{\frac{k}{2} l (\gamma)} \cdot O_k \quad \mbox{for} \; O_k \in \SO (T_x \Fh) \; .
\]
For fixed points, $\varepsilon_x = 1$ is automatic and we have:
\[
T_x \phi^t = e^{\frac{t}{2}} O_t \quad \mbox{for} \; O_t \in \SO (T_x \Fh) \; .
\]

In the number theoretical case the eigenvalues of $T_x \phi^{\log N\ep}$ on $T_x \Fh$ for $x \in \gamma_{\ep}$ would therefore be complex conjugate numbers of absolute value $N\ep^{1/2}$. If they are real then $T_x \phi^{\log N\ep}$ would simply be mutliplication by $\pm N\ep^{1/2}$. If not, the situation would be more interesting. Are the eigenvalues Weil numbers (of weight $1$)? If yes there would be some elliptic curve over $\eo_K / \ep$ involved by Tate--Honda theory. \\
{\bf 6)} It would of course be very desirable to extend the hypotheses \ref{t56} to a conjectured equality of distributions on all of $\R$. By theorem \ref{t33} we expect one contribution of the form
\[
\chi_{\Co} (\Fh \Lh , \mu) \cdot \delta_0 \; .
\]
The analogy with number theory suggests that there will also be somewhat complicated contributions from the fixed points in terms of principal values which are hard to guess at the moment. After all, even the simpler conjecture \ref{t31} has not yet been verified in the presence of fixed points!\\
{\bf 7)} If there does exist a foliated dynamical system attached to $\overline{\spec \eo_K}$ with the properties dictated by our considerations we would expect in particular that for a preferred transverse measure $\mu$ we have:
\[
\chi_{\Co} (\Fh \Lh , \mu) = - \log |d_{K / \Q}| \; .
\]
This gives some information on the space $X$ with its $\Fh \Lh$-foliation. If $K / \Q$ is ramified at some finite place i.e. if $d_{K / \Q} \neq \pm 1$ then $\chi_{\Co} (\Fh \Lh , \mu) < 0$. Now, since $\oH^2_{\Fh\Lh}$ must be one-dimensional, it follows that
\[
\chi_{\Co} (\Fh \Lh , \nu) < 0 \quad \mbox{for all non-trivial transverse measures} \; \nu \; .
\]
Hence by a result of Candel \cite{Ca} there is a Riemannian metric on $T \Fh \Lh$, such that every $\Fh \Lh$-leaf has constant curvature $-1$. Moreover $(X , \Fh\Lh)$ is isomorphic to
\[
\Oh (H , X) / \PSO (2) \; .
\]
Here $\Oh (H , X)$ is the space of conformal covering maps $u : H \to N$ as $N$ runs through the leaves of $\Fh \Lh$ with the compact open topology. See \cite{Ca} for details.

In the unramified case, $|d_{K / \Q}| = 1$ we must have $\chi_{\Co} (\Fh \Lh , \nu) = 0$ for all transverse measures by the above argument. Hence there is an $\Fh \Lh$-leaf which is either a plane, a torus or a cylinder c.f. \cite{Ca}.
\end{remarks}

\begin{punkt}
  \label{t58}
\rm In this final section we describe a simple case where the working hypothesis \ref{t56} can be proved.

Consider an unramified covering $f : M \to M$ of a compact connected orientable $d$-dimensional manifold $M$. We set
\[
\bar{M} = \lim_{\leftarrow} ( \ldots \xrightarrow{f} M \xrightarrow{f} M \to \ldots) \; .
\]
Then $\bar{M}$ is a compact topological space equipped with the shift automorphism $\bar{f}$ induced by $f$. It can be given the structure of a $C^{\infty ,\infty}$-laminated space as follows. Let $\tilde{M}$ be the universal covering of $M$. For $i \in \Z$ there exists a Galois covering
\[
p_i : \tilde{M} \longrightarrow M
\]
with Galois group $\Gamma_i$ such that $p_i = p_{i+1} \verk f$ for all $i$. Hence we have inclusions:
\[
\ldots \subset \Gamma_{i+1} \subset \Gamma_i \subset \ldots \subset \Gamma_0 =: \Gamma \cong \pi_1 (M , x_0) \; .
\]
Writing the operation of $\Gamma$ on $\tilde{X}$ from the right, we get commutative diagrams for $i \ge 0$:
\[
\begin{CD}
  \tilde{M} \times_{\Gamma} (\Gamma / \Gamma_{i+1}) @= \tilde{M} / \Gamma_{i+1} @>{\overset{p_{i+1}}{\sim}}>> M \\
@VV{\id \times \proj}V  @VV{\proj}V @VV{f}V \\
\tilde{M} \times_{\Gamma} (\Gamma / \Gamma_i) @= \tilde{M} / \Gamma_i @>{\overset{p_i}{\sim}}>> M
\end{CD}
\]
It follows that
\begin{equation}
  \label{eq:30}
  \tilde{M} \times_{\Gamma} \bar{\Gamma} \silo \bar{M}
\end{equation}
where $\bar{\Gamma}$ is the pro-finite set with $\Gamma$-operation:
\[
\bar{\Gamma} = \lim_{\leftarrow} \Gamma / \Gamma_i \; .
\]
The isomorphism (\ref{eq:30}) induces on $\bar{M}$ the structure of a $C^{\infty , \infty}$-laminated space with respect to which $\bar{f}§$ becomes leafwise smooth.

Fix a positive number $l > 0$ and let $\Lambda = l \Z \subset \R$ act on $\bar{M}$ as follows: $\lambda = l \nu$ acts by $\bar{f}^{\nu}$. Define a right action of $\Lambda$ on $\bar{M} \times \R$ by the formula
\[
(m,t) \cdot \lambda = (- \lambda \cdot m , t + \lambda) = (\bar{f}^{-\lambda / l} (m) , t + \lambda) \; .
\]
The suspension:
\[
X = \bar{M} \times_{\Lambda} \R
\]
is an $a = d+1$-dimensional $C^{\infty , \infty}$-laminated space with a one-codimensional foliation $\Fh$ as in \ref{t52}. The leaves of $\Fh$ are the fibres of the natural fibration of $X$ over the circle $\R / \Lambda$:
\[
X \longrightarrow \R / \Lambda \; .
\]
The leaves are also the images of $\bar{M} \times \{ t \}$ for $t \in \R$ under the natural projection. Translation in the $\R$-variable
\[
\phi^t [m,t'] = [m,t+t']
\]
defines an $\Fh$-compatible flow $\phi$ on $X$ which is everywhere transverse to the leaves of $\Fh$ and in particular has no fixed points.

The map
\[
\gamma \longmapsto \gamma_M = \gamma \cap (\bar{M} \times_{\Lambda} \Lambda)
\]
gives a bijection between the closed orbits $\gamma$ of the flow on $X$ and the finite orbits $\gamma_M$ of the $\bar{f}$- or $\Lambda$-action. These in turn are in bijection with the finite orbits of the original $f$-action on $M$. We have:
\[
l (\gamma) = |\gamma_M| l \; .
\]
\end{punkt}

\begin{theorem}
  \label{t59}
In the situation of \ref{t58} assume that all periodic orbits of $\phi$ are non-degenerate. Let $\Sp^n (\Theta)$ denote the set of eigenvalues with their algebraic multiplicities of the infinitesimal generator $\Theta$ of $\phi^{t*}$ on $H^n_{\Fh} (X)$. Then the trace
\[
\Tr (\phi^* \tei H^n_{\Fh} (X)) := \sum_{\lambda \in \Sp^n (\Theta)} e^{t\Theta}
\]
defines a distribution on $\R$ and the following formula holds true in $\Dh' (\R)$:
\begin{eqnarray*}
  \lefteqn{\sum^{\dim \Fh}_{n=0} (-1)^n \Tr (\phi^* \tei H^n_{\Fh} (X)) = \chi_{\Co} (\Fh\Lh , \mu) \cdot \delta_0 \; + }\\
&& \sum_{\gamma} l (\gamma) \left( \sum_{k \ge 1} \varepsilon_{\gamma} (k) \delta_{kl (\gamma)} + \sum_{k \le -1} \varepsilon_{\gamma} (|k|) \det (-T_x \phi^{kl (\gamma)} \tei T_x \Fh) \delta_{kl (\gamma)} \right) \; .
\end{eqnarray*}
Here $\gamma$ runs over the closed orbits of $\phi$ and in the sum over $k$'s any point $x \in \gamma$ can be chosen. Moreover $\chi_{\Co} (\Fh\Lh , \mu)$ is the Connes' Euler characteristic of $\Fh\Lh$ with respect to a certain canonical transverse measure $\mu$. Finally we have the formula:
\[
\chi_{\Co} (\Fh\Lh , \mu) = \chi (M) \cdot l \; .
\]
\end{theorem}

The definition of $\mu$ and the proof of the theorem will be given elsewhere.

\begin{example}
  Let $E / \F_p$ be an ordinary elliptic curve over $\F_p$ and let $\C / \Gamma$ be a lift of $E$ to a complex elliptic curve with $CM$ by the ring of integers $\eo_K$ in an imaginary quadratic field $K$. Assume that the Frobenius endomorphism of $E$ corresponds to the prime element $\pi$ in $\eo_K$. Then $\pi$ is split, $\pi \bar{\pi} = p$ and for any embedding $\Q_l \subset \C , l \neq p$ the pairs
\[
(H^*_{\et} (E \otimes \bar{\F}_p , \Q_l) \otimes \C , \Frob^*) \quad \mbox{and} \quad (H^* (\C / \Gamma , \C) , \pi^*)
\]
are isomorphic. Setting $M = \C / \Gamma , f = \pi$ we are in the situation of \ref{t58} and we find:
\[
X = (\C \times_{\Gamma} T_{\pi} \Gamma) \times_{\Lambda} \R \; .
\]
Here
\[
T_{\pi} \Gamma = \lim_{\leftarrow} \Gamma / \pi^i \Gamma \cong \Z_p
\]
is the $\pi$-adic Tate module of $\C / \Gamma$. It is isomorphic to the $p$-adic Tate module of $E$.
\end{example}

Setting $l = \log p$, so that $\Lambda = (\log p) \Z$ and passing to multiplicative time, $X$ becomes isomorphic to
\[
X \cong (\C \times_{\Gamma} T_{\pi} \Gamma) \times_{p^{\Z}} \R^*_+
\]
which may be a more natural way to write $X$. Note that $p^{\nu}$ acts on $\C \times_{\Gamma} T_{\pi} \Gamma$ by diagonal multiplication with $\pi^{\nu}$. It turns out that the right hand side of the dynamical Lefschetz trace formula established in theorem \ref{t59} equals the right hand side in the explicit formulas for $\zeta_E (s)$. Moreover the metric $g$ on $T \Fh$ given by
\[
g_{[z,y,t]} (\xi, \eta) = e^t \RRe (\xi \bar{\eta}) \quad \mbox{for} \; [z,y,t]
\; \mbox{in}\; (\C \times_{\Gamma} T_{\pi} \Gamma) \times_{\Lambda} \R 
\]
satisfies the conformality condition (\ref{eq:9}) for $\alpha = 1$. The proof of Theorem \ref{t21} can be easily adapted to $X$ above and shows that $\Theta = \halb + S$ on $\oH^1_{\Fh\Lh} (X)$ where $S$ is skew symmetric. This gives a dynamical proof for the Riemann hypotheses for $\zeta_E (s)$ along the lines that we hope for in the case of $\zeta (s)$. The construction of $(X , \phi^t)$ that we made for ordinary elliptic curves is misleading however, since it almost never happens that a variety in characteristic $p$ can be lifted to characteristic zero {\it together with its Frobenius endomorphism}. Moreover for ordinary elliptic curves the Riemann hypotheses can already be proved using Hodge cohomology of the lifted curve. This was essentially Hasse's proof.

Our present dream for the general situation is this: To an algebraic sum $\eX / \Z$ one should first attach an infinite dimensional dissipative dynamical system, possibly using $\GL_{\infty}$ in some way. The desired dynamical system should then be obtained by passing to the finite dimensional compact global attractor, c.f. \cite{La} Part I. 

\vspace{2ex}

\parbox{11.02cm}{\small
Mathematisches Institut\\
Westf. Wilhelms-Universit\"at\\
Einsteinstr. 62\\
48149 M\"unster\\
deninge@math.uni-muenster.de}
\end{document}